\documentstyle{amsppt}
\magnification=1200
\NoBlackBoxes
\def\ss{\vskip.10in}
\def\ls{\vskip.25in}
\def\P{\Bbb P}
\def\lll{$\underline{\text{\hskip.5in}}$}

\centerline{\bf Bend, Break and Count II:}
\centerline{\bf elliptics, cuspidals, linear genera}
\vskip.25in
\centerline{Z. Ran}
\centerline{Department of 
Mathematics} 
\centerline{University of California, Riverside}
\centerline{Riverside, CA  92521 USA}
\centerline{ziv\@ucrmath.ucr.edu}
\vskip.5in
In [R2] we showed how elementary considerations involving geometry on 
ruled surfaces may be used to obtain recursive enumerative formulae for 
rational plane curves.  Here we show how similar considerations may be 
used to obtain further enumerative formulae, as follow.

\proclaim{Theorem}  (i)  The number $N_d^1$ of elliptic plane curves
of degree $d$ through $3d$ general points satisfies
$$
N_d^1 = \frac{1}{12} {\binom d3}  N_d^0 + \sum_{d_1+d_2 =d}
\frac{3d_1 -2}{9} d_1 d_2 N_{d_1}^0 N_{d_2}^1 {\binom{3d-1}{3d_1-1}} ;
\tag1
$$
(ii)  The numbers $K_d^0, K_d^1$ of rational (resp. elliptic) plane curves
with 1 cusp through $3d-2$ (resp. $3d-1$) general points satisfy

$$
K_d^0 = 3N_d^0 - \sum_{d_1+d_2 = d} N_{d_1}^0 N_{d_2}^0 d_1 d_2 
[ (3d_2 -2) {\binom{3d-2}{3d_1 -2}} - \frac32 {\binom{3d-4}{3d_1-2}} ], 
\tag2
$$

$$
K_d^1 = 3N_d^1 + \frac{(d-1) (d-2)(d-4)}{8} N_d^0 + \sum (3d_1-2) d_1 d_2
{\binom{3d-1}{3d_1-1}} N_{d_1}^0 N_{d_2}^1 ;\tag3
$$
(iii)  the respective linear genera ($=$ geometric genera of generic
curve-sections) $g_d^0, g_d^1$ of the Severi varieties $V_d^0, 
V_d^1$, satisfy
$$
\align
2g_d^0 - 2 &= K_d^0 - \sum N_{d_1}^0 
N_{d_2}^0 d_1 d_2 {\binom{3d-4}{3d_1 - 2}} , \tag4
\\
2g_d^1 - 2 &= K_d^1 - \frac{9}{2} N_d^1 + \frac{(d-1) (d-2)(3d-4)}{24}
N_d^0 + \sum \frac{3d_1-2}{2} d_1 d_2 N_{d_1}^0 N_{d_2}^1 
{\binom{3d-1}{3d_1-1}} . 
\tag5
\endalign
$$
\endproclaim
Here (1) is due to 
Getzler and Pandharipande [P] while the other formulae seem new.

\subheading{Acknowledgement}  This work was begun in June `97 
under the auspices of  the `Trimestre di geometria numerativa' 
at the 
`Universit\`a di Roma Tor Vergata'.  I am grateful to the Trimestre's organizers,
Ciro Ciliberto and Corrado De Concini, as well as other participants, for 
the stimulating and congenial atmosphere, and to Rahul Pandharipande for 
telling me about his recent work.

\subheading{\S1.  Elliptics}
\ss
Let $N_d^g$ denote the Severi degree, i.e. the number of irreducible
plane curves of degree $d$ and genus $g$ through $3d+g-1$ general points.
In [R1] we gave a recursive procedure for computing $N_d^g$.  However this procedure really is for a much larger set of numbers, including $N_{d}^g$, and 
does not, as it stands, permit a recursion involving $N_d^g$ alone.  
Moreover the procedure is of exponential size in the 
number of nodes hence seems
to be at its most complicated in low-genus cases.  On the other hand for 
$N_d^0$ Kontsevich et al. have more recently given a simple recursive formula,
for which a short elementary proof is given in [R2].  
Quite recently Getzler and 
Pandharipande [P] gave an equally simple 
recursive formula for $N_d^1$, namely (1).
Their proofs are rather complicated and non-elementary.  Here we shall give an 
elementary proof in the style of [R2], which will also yield other numbers
of interest.

We begin by recalling some facts about the local geometry in 
codimension 1 of the Severi variety $V_d^g$, which may be found in [DH].
Let $\bar{B}$ denote the (closure of) the locus of curves in $V_d^g$ passing through $3d+g-2$ general points $A_1, \ldots, A_{3d+g-2} \in \P^2, B \to 
\bar{B}$ the normalization, $\bar{X} \to \bar{B}$ the tautological family 
of plane curves, $X \to \bar{X}$ the normalization and $\pi : X \to B$
the natural map.  Then $X$ is smooth and any fibre of $\pi$ is either a smooth 
irreducible curve of genus $g$, or a 1-nodal curve of arithmetic genus $g$.
The map $B \to \bar{B}$ is simply ramified over the points corresponding to 
cuspidal curves, which themselves are cusps of $\bar{B}$, and unramified 
elsewhere.  We have an exact sequence
$$
0 \to T_v \to T_X \to \pi ^*T_B \to \Cal O_{{\text sing}(\pi )} \to 0
\tag6
$$
where sing$(\pi)$ is the set of singular points of fibres and $T_v$ is the 
vertical tangent sheaf of $X$, for which (6) may be taken as a definition.
Note that $T_v$ is invertible and 
$$
c_1 (T_v) \sim - K_{X}+ (2g_d^g - 2) F,
$$
where $g_d^g = g(B), F = \pi^{-1} (pt).$

Now consider the special case $g=1$.  Then the singular fibres of $\pi$ are 
either irreducible rational curves, $\frac{(d-1)(d-2)}{2} N_d^0$ in number,
or of the form $R_i + E_i$, where $R_i$ maps to an irreducible rational curve 
of degree $d_1$ through $3d_1 -1$ of $A_1, \ldots, A_{3d-1}$ and $E_i$ maps
to an elliptic curve of degree $d_2$ through the remaining $A$'s.  Our
surface $X$ is the blow-up of a smooth 
relatively minimal elliptic surface $X'$ with 
exceptional divisor $R = \sum R_i$.  As $X$ has no multiple fibres we have,
by the well-known canonical bundle formula [BPV],
$$
K_X \sim (\omega + 2g-2) F + R
$$
where $g = g_d^1$ and $\omega$ is the degree of the $j$-function
$j : B \to M_1 = \P^1$.  Considering the fibre of $j$ over $\infty \in \P^1$,
we have
$$
\omega = \frac{1}{12} \frac{(d-1)(d-2)}{2} N_d^0. 
$$
Note also that $-\omega$ may be identified as the self-intersection of any 
section of $X'$.  In particular, if $s_{A_i} \subset X$ is the section 
corresponding (and mapping) to $A_i$, then 
$$
s_{A_i}^2 = - \omega - s_{A_i} .R . 
\tag7
$$
Now let us denote, for any function $u $ of $ d_1$,
$$
T(u) = \sum_{d_1 + d_2 = d} u (d_1) d_1 d_2 {\binom{3d-1}{3d_1 - 1}} N_{d_1}^0 
N_{d_2}^1. 
$$
Then as each $R_i$ of degree $d_1$ contains precisely $3d_1 - 1$ of the $ A$'s, we 
have $$s_{A_i}^2 = - \omega - \frac{1}{3d-1} T (3d_1 -1). $$

Now consider the ramification divisor of the natural map
$$
f: X \to \P^2.
$$
This contains a vertical part $\sum_{1}^{\kappa} F_i$ coming from the 
cuspidal curves.  Along a smooth fibre $F = \pi^{-1} (b)$, the ramification
points of $f$ are given by the residual intersection of $f(F)$ with the 
unique (mod $f(F)$) 
plane curve $C$ of degree $d$ through the nodes of $f(F)$, which
intersection is given by $A_1, \ldots, A_{3d-1}$  plus one further point 
$\sigma(b)$.  For a singular fibre $F$ the situation is similar except 
$C$ is only required to pass through the nodes of $f(F)$ not coming from 
the node of $F$; moreover, it is easy to see in case $F = E_i + R_i $ that 
the `last' point $\sigma (b)$ lies on $f(E_i)$, the elliptic part.  A moment's
thought shows that the curves $f(F)$ having a node at some $A_i$ do not 
require special consideration; what happens there is that $\bar{B}$ itself
will have a node at the point $\{f(F)\}$, with branches corresponding to the 
branches of $f(F)$ through $A_i$ and $s_{A_i} (b)$ will be the point of 
$F$ above $A_i$ lying on the branch corresponding to the branch of $\bar{B}$
containing $b$.  In any event, the points $\sigma (b)$ glue together to a 
section $\sigma$ with
$$
\multline
\sigma \sim f^* (-K_{\P^2}) + K_X - \sum_1^{3d-1} s_{A_i} - \sum_1^{\kappa} 
F_i\\
\sim 3L + (\omega + 2g-2 - \kappa) F - \sum s_{A_i} + R
\endmultline
$$
where $L = f^* {\Cal O}_{\P^2} (1)$, so that $L^2 = N_d^1$. As $\sigma \cap 
R = \emptyset$, we have
$$
\align
\sigma^2 &= - \omega \\
&= (3L + (\omega + 2g-2 - \kappa ) F - \sum s_{A_i} + R)^2 .
\endalign
$$
Simplifying, we get
$$
2g-2 - \kappa = \frac12 ((3d-9) \omega - g N_d^1 + T (3d_1 - 2 ))
\tag8
$$
On the other hand, we have (as all the $s_{A_i}$ are interchangeable)
$$
\align
(\sigma + \sum s_{A_i})^2 & =  -3d\omega - T (3d_1 -1) + 2(3d-1) \sigma. 
s_{A_1}\\
& = (3L + (\omega + 2g - 2 -\kappa ) F + R)^2 . 
\endalign
$$
Simplifying, we get
$$
\sigma .s_{A_1} = 
\frac{1}{2(3d-1)} (9N_d^1 + 9d\omega + 6d(2g-2-\kappa ) + 
T(9d_1 - 2)) . 
\tag9
$$

Now we compute $\sigma. s_{A_1}$ another way.  They key to this is the 
observation that the intersections $\sigma \cap s_{A_1}$ correspond 
precisely to those fibres $F$ on which the line bundle $\Cal O_X
(3L + R + \pi^* M)$, for any divisor $M$ on $B$, has the `same'
(i.e. isomorphic) restriction as $\Cal O_X (2s_{A_1} + s_{A_2} + 
\ldots + s_{A_{3d-1}})$; a moment's thought shows that this is valid
even for the singular fibres and that, by general choice, the intersection is 
transverse.  To make use of this observation, consider the natural 
restriction map
$$
\rho : V := \pi_* (\Cal O_X (3L+R)) \to W := \pi_* 
(\Cal O_{2s_{A_1} + s_{A_2 + \ldots + s_{A_{3d-1}}}}
(3L + R)).
$$
Then $V$ and $W$ both have rank $3d$ and $\pi (\sigma \cap s_{A_1})$
coincides with the degeneracy locus of $\rho$, hence
$$
\sigma . s_{A_1} = c_1 (W) - c_1 (V) .
\tag10
$$
Now it is easy to see that 
$$
c_1 (W) = \omega + \frac{3d+1}{3d-1} T (3d_1 - 1) . 
\tag11
$$
As for $c_1 (V)$, it may be computed from Grothendieck-Riemann-Roch [F], 
which yields 
$$
\align
c_1 (V) + 3d (1-g) &= ((1+ (3L + R) + \frac12 (3L + R)^2) \cup 
(1 - \frac{K_X}{2} + \chi (\Cal O_X) [pt])) ([X])\\
&= \frac{(3L+R)^2}{2} - \frac{K_X}{2} (3L + R) + \omega \tag12
\endalign
$$

Now (12) yields a formula for $c_1(V)$, from which the $g$ term disappears,
whence a formula for $\sigma . s_{A_1}$.  
Comparing the latter with (9) and solving for $N_d^1$, we obtain (1).

To obtain $K_d^1$, consider the natural map
$$
T_v \to f^* (T_{\P^2}) .
$$
Comparing with (6), we see that this 
map vanishes (simply) precisely at the singular points of fibres of $\pi$
plus the preimages of cusps.  Consequently
$$
\align
K_d^1 + 12 \omega + T (1) &= c_2 (f^* T_{\P^2} \otimes T_v^* ) \\
&= 3N_d^1 + 3d\omega + 3 T (d_1) - T(1) , 
\endalign
$$
which yields (3).  Finally (8) and (3) yield (5).
\ls

\subheading{\S2.  Rationals}
\ss
We now consider analogous questions for rational curves.  As above, and as 
in [R2], the idea is to study the fibred surface $\pi : X \to B$ corresponding
to the 1-parameter family of rational curves of degree $d$ through 
a set of general 
points $A_1, \ldots, A_{3d-2} \in \P^2$.  $X$ may be realized in many ways
as a blowup of a geometrically ruled surface.  One such way, $b: X \to X'= 
\P (E)$ is the blowing down of the sum $R$ of all fibre components not 
meeting $s_1 = s_{A_1}$.  Now we have
$$
s_1^2 = - m < 0, 
$$
where by [R2] we have 
$$
2m = \sum_{d_1 + d_2 = d} N_{d_1}^0 N_{d_2}^0 d_1 d_2 {\binom{3d-4}{3d_1-2}}
$$
hence we may assume $E$ has a nonvanishing section $\Cal O$ with 
$\P (\Cal O)
= s_1 \subset X'$ and $M^* = E/\Cal O$ of degree $-m$.  In fact as $X'$
possesses another section $s_{A_2}$ disjoint from $s_{A_1}$,  clearly 
$E= \Cal O \oplus M^*$, but we won't need this here.  Now let $\pi^* E^* \to 
\Cal {O}_{X'} (1)$ be the canonical invertible quotient (in the 
antiGrothendieck sense).  Then the composite $\pi^* M \to \Cal O_{X'} (1)$
vanishes precisely on $s_1$, hence
$$
c_1 (\Cal O_{X'} (1)) \sim s_1 + mf' , \quad f' = \ \text{fibre of}\ \ 
\pi'
, 
$$
hence
$$
c_1 (\Cal O_X (1)) := c_1 (b^* \Cal O_X (1)) \sim s_1 + mf \quad
f = \ \ \text{fibre of} \ \ \pi . 
$$
Consequently, letting $T_v$ as above be the vertical tangent sheaf
of $X$, we have 
$$
c_1 (T_v) = 2s_1 + mf - R .
$$
As above, the natural map
$$
T_v \to f^* T_{\P^2}
$$
vanishes precisely on singular points of fibres and cusp preimages,
hence
$$
\align
\sum N_{d_1}^0 N_{d_2}^0 d_1 d_2 {\binom{3d-2}{3d_1 -1}} + K_d^0 & = 
c_2 (f^* T_{\P^2} \otimes T_v^* ) \\ 
&= 3N_d^0 - 3dm + 3f^* \Cal O_{\P^2} (1) . R + R^2 
\endalign
$$
Simplifying, we obtain (2).  For instance, $K_3^0 = 24$, as is well known
[KS].

Finally, note that the ramification divisor  of $f: X \to \P^2$, i.e.
the degenracy locus of the differential $T_X \to f^* T_{\P^2}$, consists
now precisely of the cuspidal fibres plus the sections $s_{A_1}, \ldots,
s_{A_{3d-2}}$.  Consequently 
$$-2 s_1 -mf + R + (2g-2) f + 3f^* \Cal O (1)
\sim \sum s_{A_i} + K_d^0 f , g = g_d^0. $$  Hence taking dot product 
with $s_1$ we get $$m + 2g-2 = -m + K_d^0,$$ which yields (4).

\ls
\ls
\ls
\subheading{References}
\ss
\roster
\item"{[DH]}" S. Diaz, J. Harris:  `Geometry of Severi varieties'.
Trans. AMS. {\bf 309} (1988), 1--34. 
\ss
\item"{[F]}"  W. Fulton:  `Intersection theory'.  Springer. 
\ss
\item"{[KS]}"  S. Kleiman, R. Speiser:  `Enumerative geometry of cuspidal plane
cubics'. CMS Proc. {\bf 6} (1986), 227--268.
\ss
\item"{[P]}"  R. Pandharipande:  `Getzler's formula', preprint.
\ss
\item"{[R1]}"  Z. Ran:  `Enumerative geometry of singular plane curves'.
Invent. Math. {\bf 97}, 447--465 (1989).
\ss
\item"{[R2]}"  \lll:  `Bend, break and count', preprint (Duke eprint 9704004).

\endroster

\enddocument